\documentclass[12pt]{article}

\usepackage{amsfonts,amssymb,latexsym,setspace,subeqnarray,eqsection,epsf}

\newcommand{\proof}{\par\medskip\noindent{\sc Proof.\ }}
\newcommand{\qed}{\hfill $\Box$ \medskip \medskip}

\newtheorem{defin}{Definition}[section]

\newtheorem{prop}[defin]{Proposition}
\newtheorem{proposition}[defin]{Proposition}

\newtheorem{lemma}[defin]{Lemma}

\newtheorem{thm}[defin]{Theorem}
\newtheorem{theorem}[defin]{Theorem}

\newtheorem{corollary}[defin]{Corollary}
\newtheorem{conj}[defin]{Conjecture}
\newtheorem{conjecture}[defin]{Conjecture}

\def\reff#1{(\protect\ref{#1})}
\newcommand{\be}{\begin{equation}}
\newcommand{\ee}{\end{equation}}

\def\spose#1{\hbox to 0pt{#1\hss}}
\def\ltapprox{\mathrel{\spose{\lower 3pt\hbox{$\mathchar"218$}}
 \raise 2.0pt\hbox{$\mathchar"13C$}}}
\def\gtapprox{\mathrel{\spose{\lower 3pt\hbox{$\mathchar"218$}}
 \raise 2.0pt\hbox{$\mathchar"13E$}}}

\renewcommand{\Re}{\mathop{\rm Re}\nolimits}

\def\scrr{\mathcal{R}}

\def\R{{\mathbb R}}
\def\C{{\mathbb C}}

\begin{document}

\title{\vspace*{-2cm} On the Chromatic Roots of \\ Generalized Theta Graphs}

\author{
        \normalsize Jason Brown, Carl Hickman  \\[-1.5mm]
        \normalsize Department of Mathematics and Statistics \\[-1.5mm]
        \normalsize Dalhousie University \\[-1.5mm]
        \normalsize Halifax, Nova Scotia B3H 3J5, CANADA \\[-1.5mm]
        \normalsize {\tt brown@mathstat.dal.ca},
                    {\tt chickman@mathstat.dal.ca} \\
        ~ \\
        \normalsize Alan D.~Sokal \\[-1.5mm]
        \normalsize Department of Physics \\[-1.5mm]
        \normalsize New York University \\[-1.5mm]
        \normalsize New York, NY 10003, USA \\[-1.5mm]
        \normalsize {\tt sokal@nyu.edu} \\
        ~ \\
        \normalsize David G. Wagner \\[-1.5mm]
        \normalsize Department of Combinatorics and Optimization \\[-1.5mm]
        \normalsize University of Waterloo \\[-1.5mm]
        \normalsize Waterloo, Ontario N2L 3G1, CANADA \\[-1.5mm]
        \normalsize {\tt dgwagner@math.uwaterloo.ca}
        ~\\[5mm]
}

\date{October 31, 2000}

\maketitle
\thispagestyle{empty}   

\begin{abstract}
The generalized theta graph $\Theta_{s_1,\ldots,s_k}$
consists of a pair of endvertices joined by
$k$ internally disjoint paths of lengths $s_1,\ldots,s_k \ge 1$.
We prove that the roots of the chromatic polynomial
$\pi(\Theta_{s_1,\ldots,s_k},z)$ of a $k$-ary generalized theta graph
all lie in the disc $|z-1| \le [1 + o(1)] \, k/\log k$,
uniformly in the path lengths $s_i$.
Moreover, we prove that $\Theta_{2,\ldots,2} \simeq K_{2,k}$
indeed has a chromatic root of modulus $[1 + o(1)] \, k/\log k$.
Finally, for $k \le 8$ we prove that the generalized theta graph
with a chromatic root that maximizes $|z-1|$
is the one with all path lengths equal to 2;
we conjecture that this holds for all $k$.
\end{abstract}

\vspace{0.5cm}

\clearpage

\section{Introduction}

The {\em chromatic polynomial}\/ $\pi(G,z)$ of a
(finite undirected) graph $G=(V,E)$
is the number of proper $z$-colorings of the vertices of $G$,
i.e., of functions $f\colon\, V\rightarrow \{1,2,\ldots ,z\}$
such that $uv\in E$ implies $f(u)\neq f(v)$.
As chromatic polynomials are indeed polynomials \cite{Read_68,Read_88},
their roots (called {\em chromatic roots}\/)
have been intensively studied by combinatorialists \cite{Read_88}.
Chromatic polynomials also arise in statistical physics as
zero-temperature limits of the partition function of the
$z$-state Potts antiferromagnet on $G$ (see e.g.\ \cite{Sokal_le70});
and since the complex zeros of the partition function
are closely linked to phase transitions \cite{Yang-Lee_52},
physicists (or at least a few of them)
have likewise been much interested in locating chromatic roots
\cite{Shrock_99g,Salas-Sokal_transfer1}.

Very recently, one of us proved the following,
which confirms a 1972 conjecture of Biggs, Damerell, and Sands \cite{biggs}:

\begin{thm}[Sokal \cite{Sokal_00a}]
 \label{thm1.1}
If $G$ is a graph of maximum degree $k$,
then every chromatic root $z$ of $G$ lies in the disc $|z|< 7.963907\,k$.
\end{thm}

\noindent
While the constant $7.963907$ found in \cite{Sokal_00a} can likely be improved,
the linearity of the bound is best possible,
since the complete graph $K_{k+1}$ (which has maximum degree $k$)
has a chromatic root at $z=k$.\footnote{
   Surprisingly, the complete graph $K_{k+1}$ is {\em not}\/
   the extremal graph for this problem.
   A nonrigorous (but probably rigorizable) asymptotic analysis,
   confirmed by numerical calculations, shows \cite{Salas-Sokal_in_prep}
   that the complete bipartite graph $K_{k,k}$
   has a chromatic root $\alpha k + o(k)$, where
   $\alpha = - 2 / W(-2/e)  \approx 0.678345 + 1.447937 i$;
   here $W$ denotes the principal branch of the Lambert $W$ function
   (the inverse function of $w \mapsto w e^w$) \cite{Corless_96}.
   So the constant in Theorem~\ref{thm1.1} cannot be better than
   $|\alpha| \approx 1.598960$.
   Gordon Royle (private communication) has conjectured that
   among all graphs of maximum degree $k$ ($k \ge 4$),
   the graph with the largest modulus of a chromatic root is $K_{k,k}$.
   We conjecture that this holds also when $|z|$ is replaced by $|z-1|$,
   as is done in the present paper.
}
Nevertheless, suitably restricted {\em subclasses}\/ of graphs
might well satisfy a sublinear bound.
For example, we conjecture that a sublinear bound in terms of
maximum degree holds for all {\em series-parallel}\/ graphs,
and conceivably even for all {\em planar}\/ graphs.

A complementary approach is to bound the chromatic roots
in terms of the {\em corank}\/ (or {\em cyclomatic number}\/)
of the graph.  Recall that a connected graph with
$n$ vertices and $m$ edges has corank $m-n+1$.
One of us has recently proven:

\begin{thm}[Brown \cite{brown1}]
 \label{thm1.2}
If $G$ is a graph of corank $k\geq1$,
then every chromatic root $z$ of $G$ lies in the disc $|z-1|\leq k$.
\end{thm}

\noindent
Unlike the bound in terms of maximum degree, however,
it is not known whether linear growth with corank is best possible.
Indeed, we suspect that it is not (see Section~\ref{sec_conclusion}).

Note also that Theorem~\ref{thm1.1} is in some sense radically
``more powerful'' than Theorem~\ref{thm1.2}, in that maximum degree
is a ``local'' quantity (i.e.\ a max over vertices),
while corank is a ``global'' quantity (i.e.\ a sum over vertices).
Thus, if $G$ is (for example) a regular graph of degree $\Delta \ge 3$
with $n$ vertices, then the maximum degree is always $\Delta$
(independent of $n$), while the corank is $n(\Delta-2)/2 + 1$
(hence grows linearly with $n$).
More generally, for any 2-connected loopless graph we have
${\rm cr}(G) \ge \Delta(G) - 1$, where ${\rm cr}(G)$ is the corank
and $\Delta(G)$ is the maximum degree.\footnote{
   {\sc Proof.}  Let $G$ be a 2-connected loopless (multi)graph
   with $n \ge 1$ vertices and $m$ edges.  The proof is by induction on $m$.
   If $m=0$, then $G=K_1$ and ${\rm cr}(G) = \Delta(G) = 0$.
   If $m=1$, then $G=K_2$ and ${\rm cr}(G) = \Delta(G) -1 = 0$.
   So suppose $m \ge 2$.
   Then, by 2-connectedness, $G$ cannot have any degree-1 vertices.
   If $G$ has a degree-2 vertex where the two edges connect to
   distinct vertices, then by contracting that vertex we obtain a
   loopless graph $G'$ that is homeomorphic to $G$ (hence 2-connected)
   and satisfies $m(G') = m(G) -1$, ${\rm cr}(G') = {\rm cr}(G)$ and
   $\Delta(G') = \Delta(G)$;  so we can apply the inductive hypothesis.
   Likewise, if $G$ has somewhere $k \ge 2$ parallel edges between the
   same pair of distinct vertices, then by replacing those parallel edges
   by a single edge we obtain a graph $G'$ that is 2-connected and satisfies
   $m(G') = m(G) - (k-1)$, ${\rm cr}(G') = {\rm cr}(G) - (k-1)$
   and $\Delta(G') \ge \Delta(G) - (k-1)$;
   so again we can apply the inductive hypothesis.
   Finally, if $G$ is a simple graph of minimum degree $\ge 3$
   and maximum degree $\Delta$,
   then $\Delta \le n-1$ and $m \ge [\Delta + 3(n-1)]/2$;
   this implies ${\rm cr}(G) = m-n+1 \ge [\Delta + (n-1)]/2 \ge \Delta$,
   which is even stronger than what we want to prove.

   The same proof shows that if $G$ is 2-connected and not series-parallel,
   then ${\rm cr}(G) \ge \Delta(G)$:
   it suffices to observe that the two reductions
   (contraction of a degree-2 vertex and replacement of parallel edges
    by a single edge)
   preserve the property of being non-series-parallel,
   and that $K_2$ is series-parallel;
   therefore, after all the reductions we must be left with
   a simple graph of minimum degree $\ge 3$.
   [Remark:  The smallest non-series-parallel graph, $K_4$,
    has ${\rm cr}(K_4) = \Delta(K_4) = 3$.]
}
Nevertheless, for some graphs Theorem~\ref{thm1.2} will give
a sharper bound than Theorem~\ref{thm1.1},
because of its prefactor 1 in place of 7.963907.

\begin{figure}[t]
\begin{center}
\epsfxsize = 0.5\textwidth
\leavevmode\epsffile{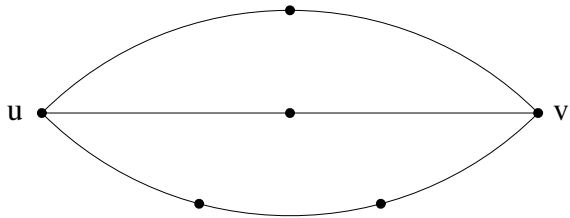}
\end{center}
\vspace*{0cm}
\caption{The graph $\Theta_{2,2,3}$.}
\end{figure}

In this paper, we shall examine a subclass of series-parallel graphs
for which we are able to prove a sublinear bound
on the chromatic roots in terms of {\em both}\/ corank and maximum degree.
The graphs in question are the {\em generalized theta graphs}\/
$\Theta_{s_1,\ldots,s_k}$,
which are formed by taking a pair of vertices $u,v$
(called the {\em endvertices}\/)
and joining them by $k$ internally disjoint paths
of lengths $s_1,\ldots,s_k \ge 1$.
(A generalized theta graph with three paths is traditionally called a
 {\em theta graph}\/ without adjectives.
 However, the rationale for singling out the case $k=3$
 seems more philological than mathematical,
 so we shall often drop the adjective ``generalized''
 when referring to $k$-ary theta graphs.)
For brevity, we denote by $\Theta^{(s,k)}$
the $k$-ary theta graph whose path lengths are all equal to $s$.
A $k$-ary theta graph clearly has maximum degree $k$ and corank $k-1$
(except for the trivial case $k=1$ with $s_1 \ge 2$,
 which has maximum degree 2).
The main result of the present paper is the following:

\begin{thm}
 \label{thm1.3}
The chromatic roots of any $k$-ary generalized theta graph
lie in the disc $|z-1| \le [1 + o(1)] \, k/\log k$,
where $o(1)$ denotes a constant $C(k)$ that tends to zero as $k\to\infty$.
(The precise bound is given in Theorem~\ref{thm3.2} and
 Proposition~\ref{prop_R2k}.)
\end{thm}

\noindent
Moreover, we shall prove that this bound is asymptotically saturated
by the graph $\Theta^{(2,k)}$ with all path lengths equal to 2
(which is isomorphic to the complete bipartite graph $K_{2,k}$).
Finally, we shall prove that for $k \le 8$,
the $k$-ary generalized theta graph that maximizes $|z-1|$
over its chromatic roots $z$ is the one with all path lengths equal to 2;
we conjecture that this holds for all $k$.

More generally, for any finite graph $G$ let us define
\be
   \rho(G)   \;=\;  \max\{ |z-1| \colon\; \pi(G,z) = 0 \}
   \;,
\ee
and let us write specifically
\begin{eqnarray}
   \rho(s_1,\ldots,s_k)   & = &   \rho(\Theta_{s_1,\ldots,s_k})   \\
   \rho^{(s,k)}   & = &   \rho(\Theta^{(s,k)})
\end{eqnarray}
Broadly speaking, this paper will be concerned with finding
upper and lower bounds on $\rho(s_1,\ldots,s_k)$.


This paper is concerned, therefore, with the behavior of the chromatic roots
of $\Theta_{s_1,\ldots,s_k}$ when $k$ is {\em fixed}\/ (albeit large).
A very different situation arises when $k$ is unbounded:
indeed, one of us has shown elsewhere \cite{Sokal_00b}
that the chromatic roots of the graphs $\Theta^{(s,k)}$ ($s,k \ge 2$),
taken together, are dense in the whole complex plane
except possibly for the disc $|z-1| < 1$!

The plan of this paper is as follows:
In Section~\ref{computation} we compute the chromatic polynomial
of $\Theta_{s_1,\ldots,s_k}$.
In Section~\ref{bound} we prove Theorem~\ref{thm1.3}.
In Section~\ref{smalltheta} we use a different method to
obtain sharper (at least for small $k$) upper bounds
on the chromatic roots of $\Theta_{s_1,\ldots,s_k}$;
as a corollary, we show that for $k \le 8$ the graph that maximizes
$\rho(s_1,\ldots,s_k)$ is the one with all path lengths equal to 2.
In Section~\ref{K2k} we prove that $\Theta^{(2,k)} \simeq K_{2,k}$
indeed has a chromatic root of magnitude $[1 + o(1)] \, k/\log k$,
so that Theorem~\ref{thm1.3} is asymptotically sharp.
Finally, in Section~\ref{sec_conclusion} we discuss our results
and make some conjectures.

\section{Chromatic Polynomial of Generalized Theta Graphs} \label{computation}

It is well known that the chromatic polynomial of a cycle of length $n$
is given by
\begin{eqnarray}
\pi(C_{n},z) & = & (z-1)^{n} + (-1)^{n}(z-1)   \;.   \label{cycle}
\end{eqnarray}
Recall also the well-known {\em addition-contraction formula}\/
\begin{eqnarray}
\pi(G,z) & = & \pi(G + e,z) + \pi(G \bullet e,z) \label{delCon}
\end{eqnarray}
where $e$ is any edge not in $G$ and
$G \bullet e$ denotes the graph obtained from $G$ by contracting
the endpoints of $e$.
Finally, if two graphs $G$ and $H$ overlap in a complete graph on $l$ vertices,
then
\begin{eqnarray}
\pi(G \cup H, z) & = & \frac{\pi(G,z) \, \pi(H,z)}{z(z-1) \cdots (z-l+1)}
    \;.  \label{union}
\end{eqnarray}

Observe now that if one adjoins to a generalized theta graph
$\Theta_{s_1,\ldots,s_k}$ a new edge $e$ between the two endvertices,
the resulting graph is a collection of cycles (of lengths $s_i +1$)
that overlap in a complete graph $K_2$ (namely the edge $e$),
so that
\be
   \pi(\Theta_{s_1,\ldots,s_k} + e, z)   \;=\;
   \frac{\prod\limits_{i=1}^{k}{ [(z-1)^{s_i +1}+(-1)^{s_i +1}(z-1)]}}
        {[z(z-1)]^{k-1}}
   \;.
\ee
Likewise, contraction of the endvertices of $\Theta_{s_1,\ldots,s_k}$
yields a collection of cycles (of lengths $s_i$)
that overlap in a complete graph of order 1,
so that
\be
   \pi(\Theta_{s_1,\ldots,s_k} \bullet e, z)   \;=\;
   \frac{\prod\limits_{i=1}^{k}{ [(z-1)^{s_i}+(-1)^{s_i}(z-1)]}}
        {z^{k-1}} 
   \;.
\ee
It follows that the chromatic polynomial of
$\Theta_{s_{1},\ldots ,s_{k}}$ is given by
\begin{subeqnarray}
\pi(\Theta_{s_{1},\ldots , s_{k}},z) & = &
  \frac{\prod\limits_{i=1}^{k}{ [(z-1)^{s_i +1}+(-1)^{s_i +1}(z-1)]}}
       {[z(z-1)]^{k-1}}   \;+\;  \nonumber \\
 & & \qquad \frac{\prod\limits_{i=1}^{k}{ [(z-1)^{s_i}+(-1)^{s_i}(z-1)]}}
                 {z^{k-1}}  \\[2mm]
 & = & \frac{(-1)^{\sum\limits_{i=1}^{k}{s_i}-1}(1-z)}{z^{k-1}}
       \left[ \prod_{i=1}^{k}{[(1-z)^{s_i}-1]} 
              \;-\; \right. \nonumber \\
 & & \qquad\qquad \left. (1-z)^{k-1}\prod_{i=1}^{k}{[(1-z)^{s_{i}-1}-1]} \right]
     \;. \qquad
 \label{theta}
\end{subeqnarray}
Therefore, we need only concern ourselves with the roots of
\begin{equation}
 f_{s_{1},\ldots ,s_{k}}(y)  \;=\;   \prod_{i=1}^{k}(y^{s_i}-1) -
y^{-1} \prod_{i=1}^{k} (y^{s_i}-y)
   \;,
 \label{thetaroots}
\end{equation}
where $y=1-z$.
All our subsequent calculations will be expressed
in terms of the variable $y$.

Let us now dispose of some trivial cases.
If $k=1$, the theta graph $\Theta_{s_1}$
is isomorphic to the path $P_{s_1}$,
so that its chromatic roots are 0 and 1.
If $k=2$, the theta graph $\Theta_{s_1,s_2}$
is isomorphic to the cycle $C_{s_1+s_2}$,
so that its chromatic roots (other than $z=1$)
all lie on the circle $|z-1| = 1$.
We shall therefore assume henceforth that $k \ge 3$.

If one or more of the path lengths $s_i$ equals 1,
then the second product in \reff{thetaroots} vanishes,
and all the chromatic roots (other than $z=1$)
again lie on the circle $|y| = 1$, i.e.\ $|z-1| = 1$.
This can also be seen immediately from \reff{cycle} and \reff{union},
as in this case the graph $\Theta_{s_1,\ldots,s_k}$ is a
$K_2$--bond of cycles (i.e.\ a collection of cycles overlapping
in a single edge).
We shall therefore assume henceforth that all path lengths $s_i$ are $\ge 2$.

\bigskip

\noindent
{\bf Remarks.}
1. This method for computing $\pi(\Theta_{s_1,\ldots,s_k}, z)$
was employed previously by Read and Tutte \cite[pp.~29--30]{Read_88}
for the case $k=3$.

2. See \cite[Section 2]{Sokal_00b} for a systematic method for computing
the Potts-model partition function $Z_G(z, \{v_e\})$
[which generalizes the chromatic polynomial $\pi(G,z)$]
for any series-parallel graph,
along with an explicit formula for generalized theta graphs.

\section{Bounding the Chromatic Roots of Generalized Theta Graphs}
    \label{bound}

We proceed now to prove Theorem~\ref{thm1.3},
i.e.\ to show that the roots of $f_{s_1,\ldots,s_k}(y)$
are bounded in modulus by $[1 + o(1)] \, k/\log k$.
As just noted, it suffices to consider the case $k \ge 3$
and $s_1, \ldots, s_k \ge 2$.

Obviously $y=1$ (corresponding to $z=0$) is a root of $f_{s_1,\ldots,s_k}(y)$.
All the other roots satisfy
\be
   \prod_{i=1}^k  {y^{s_i} - y  \over  y^{s_i} - 1}
   \;=\;
   y \;.
 \label{eq3.1}
\ee
Therefore, to show that a given number $y \in \C$ is not a root,
it suffices to show that the left-hand side of \reff{eq3.1}
is strictly smaller in modulus than the right-hand side.
We shall do this in the crudest possible way,
taking account only of the magnitude of $y$,
i.e.\ throwing away all phase information.
So, define
\be
   X_s(R)  \;=\;
   \sup\limits_{|y| = R}  \, \left| {y^s - y  \over  y^s - 1} \right|
 \label{def_Xs}
\ee
for integer $s \ge 2$ and real $R > 1$.

\begin{lemma}
   \label{lemma3.1}
\quad\par
\begin{itemize}
   \item[(a)]
For integer $s \ge 2$ and real $R > 1$, we have
\be
   X_s(R)  \;\le\;  \widetilde{X}_s(R) \,\equiv\, {R^s + R  \over  R^s - 1}
   \;,
 \label{bound_Xs}
\ee
with equality for $s$ even
[attained uniquely in \reff{def_Xs} when $y = -R$]
and strict inequality for $s$ odd.
   \item[(b)]
For integer $s \ge 2$, both $X_s(R)$ and $\widetilde{X}_s(R)$
are strictly decreasing functions of $R$ on the interval $1 < R < \infty$.
They tend to $+\infty$ as $R \downarrow 1$, and to 1 as $R \to \infty$.
   \item[(c)]
For $R > 1$, $\widetilde{X}_s(R)$ is a strictly decreasing function of $s$
on $s \ge 2$.
\end{itemize}
\end{lemma}

\proof
(a)  The inequality \reff{bound_Xs} is trivial.
Equality holds if and only if we simultaneously have
$y^{s-1}$ negative real and $y^s$ positive real.
Writing $y=Re^{i\theta}$ with $0 \le \theta < 2\pi$,
a simple computation shows that this occurs if and only if
$s$ is even and $\theta = \pi$.
Thus, for $s$ odd we have $|(y^s - y)/(y^s - 1)| < (R^s + R)/(R^s - 1)$
for all $y$ on the circle of radius $R$;
and since this circle is compact, the maximum is attained;
so \reff{bound_Xs} is a strict inequality.

(b,c)  The function $g(y) = (y^s - y)/(y^s - 1)$ is analytic
in the domain $|y| > 1$, including at the point $y=\infty$;
so by the maximum modulus theorem, its maximum on the domain $|y| \ge R$
(with $R > 1$) is attained only on the boundary $|y| = R$.
This proves (invoking again the compactness of circles)
that $X_s(R)$ is a strictly decreasing function of $R$ for $R > 1$.
Finally, simple calculus shows that $\widetilde{X}_s(R) = 1 + (R+1)/(R^s -1)$
is a strictly decreasing function of $R$ on $R>1$;
and it is trivially a strictly decreasing function of $s$.
\qed

\medskip
\noindent
{\bf Remarks.}
1.  For $s$ odd, the maximum in \reff{def_Xs} lies at
a nontrivial angle $\theta$,
which tends to $\pi \pm \pi/(s-1)$ as $R \to \infty$
and to $\pi \pm \pi/s$ as $R \downarrow 1$.\footnote{
   {\sc Proof.}
   Let $y = R e^{i\theta}$.
   As $R \to \infty$ we have
   $$
      {y^s - y  \over  y^s - 1}   \;=\;
      1 \,-\, {1 \over y^{s-1}} \,+\, {1 \over y^s} \,+\,
               O\!\left( {1 \over y^{2s-1}} \right)
      \;.
   $$
   The sum of the first two terms has maximum modulus when
   $e^{i\theta}$ is an $(s-1)^{st}$ root of $-1$,
   i.e.\ when $\theta = (2k+1)\pi/(s-1)$ with $0 \le k \le s-2$;
   and among these, the third term has maximum real part when
   $k = (s-1)/2$ and $(s-3)/2$,
   i.e.\ when $\theta = \pi \pm \pi/(s-1)$.
   As $R \downarrow 1$, the denominator dominates,
   so we want $e^{i\theta}$ to be an $s^{th}$ root of unity,
   i.e.\ $\theta = 2\pi k/s$ with $0 \le k \le s-1$;
   and among these, the numerator has maximum modulus when
   $y^{s-1}$ has the largest negative real part,
   which occurs when $k = (s \pm 1)/2$,
   i.e.\ when $\theta = \pi \pm \pi/s$.
}

2.  For $s$ odd, both sides of \reff{bound_Xs} have the same leading
asymptotic behavior as $R \to \infty$, namely $1 + 1/R^{s-1} + O(1/R^s)$.
So the inequality \reff{bound_Xs} is asymptotically sharp.

\bigskip


\bigskip

It follows that if
\be
   \prod_{i=1}^k  \widetilde{X}_{s_i}(R) \;\le\; R
   \;,
 \label{prod_Xtilde}
\ee
then $f_{s_1,\ldots,s_k}(y)$ has no roots in the region $|y| > R$.
Now, the left-hand side of \reff{prod_Xtilde} is a strictly decreasing
function of $R$ on $1 < R < \infty$,
which tends to $+\infty$ as $R \downarrow 1$ and to 1 as $R \to\infty$;
while the right-hand side is a strictly increasing function of $R$,
which tends to 1 as $R \downarrow 1$ and to $+\infty$ as $R \to\infty$.
So there is a unique $R = \scrr(s_1,\ldots,s_k)$ where they are equal.
Moreover, since the left-hand side of \reff{prod_Xtilde}
is a strictly decreasing function of each $s_i$,
it follows that $\scrr(s_1,\ldots,s_k)$ is a strictly decreasing function
of each $s_i$;  in particular,
$\scrr(s_1,\ldots,s_k) \le \scrr(2,\ldots,2) \equiv \scrr^{(2,k)}$.
[When all the $s_i$ are equal, we write
 $\scrr(s,\ldots,s) = \scrr^{(s,k)}$.]
In summary, we have proven:

\begin{theorem}
   \label{thm3.2}
All the chromatic roots of $\Theta_{s_1,\ldots,s_k}$ lie in the disc
$|z-1| \le \scrr(s_1,\ldots,s_k)$.
In particular, they all lie in the disc $|z-1| \le \scrr^{(2,k)}$.
\end{theorem}

To complete the proof of Theorem~\ref{thm1.3},
we need to obtain an upper bound on $\scrr^{(2,k)}$.
Since $X_2(R) = \widetilde{X}_2(R) = R/(R-1)$,
we see that $\scrr^{(2,k)}$ is the unique root in $(1,\infty)$ of the equation
\be
   \left( {R \over R-1} \right) ^{\! k}   \;=\;   R
 \label{eq3_Ra}
\ee
or equivalently
\be
   \left( {R \over R-1} \right) ^{\! k-1}   \;=\;   R-1
 \label{eq3_Rb}
\ee
or
\be
   R^{k-1}   \;=\;  (R-1)^k
   \;.
 \label{eq3_Rc}
\ee
The asymptotic behavior as $k \to\infty$
of the solution to \reff{eq3_Ra}--\reff{eq3_Rc}
is surprisingly intricate
and involves the Lambert $W$ function \cite{Corless_96},
i.e.\ the inverse function to $w \mapsto w e^w$.
Here we shall limit ourselves to proving some elementary inequalities,
without deriving the full asymptotics.
For further details, see Section~\ref{K2k} below and \cite{Sokal_trinomials}.

\begin{lemma}
   \label{lemma_lambert}
For $x \ge 0$, let $W(x)$ be the unique real solution of $W(x) e^{W(x)} = x$.
Then:
\begin{itemize}
   \item[(a)] $W(x)$ is a continuous (in fact, real-analytic)
      and strictly increasing function of $x$,
      satisfying $W(0) = 0$ and $\lim\limits_{x\to +\infty} W(x) = +\infty$.
   \item[(b)] $x/W(x) = e^{W(x)}$ is a continuous (in fact, real-analytic)
      and strictly increasing function of $x$,
      satisfying $\lim\limits_{x\downarrow 0} x/W(x) = 1$ and
       $\lim\limits_{x\to +\infty} x/W(x) = +\infty$.
   \item[(c)] For all $x > e$,
\be
   \log x \,-\, \log\log x  \;<\;  W(x)  \;<\;  \log x   \;.
\ee
\end{itemize}
\end{lemma}

\proof
%
(a) follows immediately from the fact that $w e^w$ is a real-analytic
function of $w \ge 0$ with everywhere strictly positive derivative,
which runs from 0 to $+\infty$ as $w$ runs from 0 to $+\infty$.
(b) follows from (a) and the defining equation $W(x) e^{W(x)} = x$.
To prove (c), note that $W(e) = 1$, so that for $x > e$ we have $W(x) > 1$
by (a);  it follows that $W(x) = \log x - \log W(x) < \log x$
and hence also that $W(x) = \log x - \log W(x) > \log x - \log\log x$.
\qed

\begin{proposition}
   \label{prop_R2k}
For $k \ge 2$ we have
\be
   {k \over W(k)}  \;<\;
   \scrr^{(2,k)}  \;<\;
   {k-1 \over W(k-1)}  \,+\, 1
   \;.
 \label{eq_prop_R2k_a}
\ee
In particular, for $k \ge 3$ we have
\be
   {k \over \log k}  \;<\;
   \scrr^{(2,k)}  \;<\;
   {k-1 \over \log (k-1) \,-\, \log\log (k-1)}  \,+\, 1
   \;.
 \label{eq_prop_R2k_b}
\ee
\end{proposition}

\bigskip
\noindent
It suffices to prove \reff{eq_prop_R2k_a},
as \reff{eq_prop_R2k_b} then follows from Lemma~\ref{lemma_lambert}.

\par\medskip\noindent
{\sc Proof of lower bound.}
Since the left (resp.\ right) side of \reff{eq3_Ra} is a strictly decreasing
(resp.\ strictly increasing) function of $R$,
it suffices to prove that the left side is larger than the right side
when $R = k /W(k)$.
Working with reciprocals, we have
\be
   \left( {R-1 \over R} \right) ^{\! k}
   \;=\;
   \left( 1 \,-\, {W(k) \over k} \right) ^{\! k}
   \;<\; \exp[-W(k)]
   \;=\;
   {W(k) \over k}
   \;=\;
   {1 \over R}
   \;.
\ee
(In the inequality step we used $0 < W(k)/k < 1$.)

\par\medskip\noindent
{\sc Proof of upper bound.}
Since the left (resp.\ right) side of \reff{eq3_Rb} is a strictly decreasing
(resp.\ strictly increasing) function of $R$,
it suffices to prove that the left side is smaller than the right side
when $R = (k-1)/W(k-1) \,+\, 1$:
\be
   \left( {R \over R-1} \right) ^{\! k-1}
   \;=\;
   \left( 1 \,+\, {W(k-1) \over k-1} \right) ^{\! k-1}
   \;<\; \exp[W(k-1)]
   \;=\;
   {k-1 \over W(k-1)}
   \;=\;
   R-1
   \;.
\ee
\qed

{\bf Remarks.}
1.  The precise large-$x$ asymptotics of $W(x)$ is
\cite{deBruijn_61,Comtet_70,Corless_96}
\begin{subeqnarray}
   W(x)  & = &   \log x \,-\, \log\log x \,+\,
                 \sum_{n=1}^\infty \sum_{k=1}^n
                      (-1)^{n+1} \, {s(n, n-k+1) \over k!} \,
                      {(\log\log x)^k \over (\log x)^n}
      \\[2mm]
   & = &  \log x \,-\, \log\log x \,+\,  {\log\log x \over \log x}
          \,+\, {(\log\log x)^2 \over 2 (\log x)^2}
          \,-\, {\log\log x \over (\log x)^2}
          \,+\, O\!\left( {(\log\log x)^3 \over (\log x)^3} \right)
     \nonumber \\
 \label{lambert_asymptotics}
%
%
\end{subeqnarray}
where the $s(n,m)$ are the Stirling numbers of the first kind.
Surprisingly, this series is {\em convergent}\/ for sufficiently large $x$
\cite{deBruijn_61,Corless_96}.

2.  The precise large-$k$ asymptotics of $\scrr^{(2,k)}$ is
\cite{Sokal_trinomials}
\be
   \scrr^{(2,k)}  \;=\;
   {k \over W(k)} \,+\, {W(k)  \over  2 [1 + W(k)]}  \,+\,
   O\!\left( {W(k) \over k} \right)
 \;.
\ee
Again, this is a {\em convergent}\/ series for sufficiently large $k$;
it can be obtained from equations
\reff{def_zeta}/\reff{def_xi}/\reff{series2_xi_v}/\reff{g1_v} below.

\section{Generalized Theta Graphs with $\le 8$ Paths} \label{smalltheta}

Numerical computations suggest that among all $k$-ary theta graphs
$\Theta_{s_1,\ldots , s_k}$,
the one with a chromatic root that maximizes $|z-1|$ is the
graph with all path lengths $s_i$ equal to 2,
i.e. the graph $\Theta^{(2,k)} \simeq K_{2,k}$.
The general bounds of the previous section are not, however,
strong enough to prove this conjecture.
Nevertheless, by different techniques we shall show the validity
of this conjecture for all $k \leq 8$.
As before, it suffices to consider $k \ge 3$ and $s_1,\ldots,s_k \ge 2$.
We recall that $\rho(s_{1},\ldots ,s_{k})$ denotes
the maximum modulus of a root of $f_{s_{1},\ldots, s_{k}}(y)$.

Our method is based on the following trivial bound
(see e.g.\ \cite[Theorem 27.1]{marden}):

\begin{proposition}
   \label{prop4.1}
Let $P(y) = \sum\limits_{j=0}^n a_j y^j$
be a polynomial of degree $n$ (so that $a_n \neq 0$),
and let $R$ be the unique nonnegative real solution of
\be
   |a_n| R^n \,-\, \sum\limits_{j=0}^{n-1} |a_j| R^j   \;=\;   0   \;.
 \label{eq4.1}
\ee
Then all the roots of $P$ lie in the disc $|y| \le R$.
Moreover, in \reff{eq4.1} the numbers $|a_j|$ for $0 \le j \le n-1$
can be replaced by any numbers $b_j \ge |a_j|$;
this only makes the bound weaker.
\end{proposition}

\noindent
At first sight it is surprising that such a crude estimation method
--- which throws away all the sign or phase information
in the coefficients of $P$ ---
could yield reasonably sharp results.
And indeed, we do not entirely understand why it works so well
in our application --- but it does!
Here is the key trick:
since the polynomial $f_{s_1,\ldots,s_k}(y)$
defined in (\ref{thetaroots}) is divisible by $(y-1)^k$,
we are free to pull out a factor $(y-1)^l$ with any $0 \le l \le k$
before applying Proposition~\ref{prop4.1}.
It turns out that the right choice is to take $l=1$
(in the Remark below we shall give some intuition as to why this is
a good choice).
That is, we define the polynomial $\phi_{s_{1},\ldots,s_{k}}(y)$ by
\begin{equation}
  \phi_{s_{1},\ldots ,s_{k}}(y)  \;=\;
  {f_{s_{1},\ldots ,s_{k}}(y)  \over  y-1}
  \;.
 \label{phipoly}
\end{equation}
Let $[k]$ denote the set $\{1,\ldots,k\}$,
and let ${{[k]} \choose {l}}$ denote the set of all subsets of $[k]$
of cardinality $l$.
For a subset $X = \{ i_{1},i_{2},\ldots,i_{l} \}$ of $[k]$, let us define
\[ s_{X} = \sum_{j=1}^{l} s_{i_{j}}.\]
We then have:


\begin{subeqnarray}
   f_{s_1,\ldots,s_k}(y)  & = &
      \prod_{i=1}^{k}(y^{s_i}-1) - y^{-1} \prod_{i=1}^{k}(y^{s_i}-y)
           \\[2mm]
 &=&  \sum_{m=0}^k (-1)^m \sum_{X \subseteq {[k] \choose {k-m}}} y^{s_X}
        \,-\,
      \sum_{m=0}^k (-1)^m \sum_{X \subseteq {[k] \choose {k-m}}} y^{s_X +m-1}
           \\[2mm]
 &=&  \sum_{m=0}^k (-1)^m \sum_{X \subseteq {[k] \choose {k-m}}}
           y^{s_X} (1 - y^{m-1})
           \\[2mm]
 &=&  y^{(\sum{s_{i}}) - 1} (y-1) \,-\,
           \sum_{m=2}^k (-1)^m \sum_{X \subseteq {[k] \choose {k-m}}}
           y^{s_X} (y^{m-1} - 1)
\end{subeqnarray}
and hence

\be
   \phi_{s_1,\ldots,s_k}(y)   \;=\;
      y^{(\sum{s_{i}}) - 1} \,-\,
      \sum_{m=2}^k (-1)^m \sum_{X \subseteq {[k] \choose {k-m}}}
           y^{s_X} (1 + y + \cdots + y^{m-2})
      \;.
 \label{phiform}
\ee


\bigskip
\noindent
{\bf Remark.}
Note that in $\phi_{s_1,\ldots,s_k}(y)$,
the first subleading term with a $+$ sign (which comes from $m=3$)
is down by a factor $y^{s_1 + s_2 + s_3 - 2} = y^{\ge 4}$
compared to the leading term,
where $s_1, s_2, s_3 \ge 2$ are the three smallest path lengths.
Since $y > 1$, this helps Proposition~\ref{prop4.1} to be close to sharp.

\bigskip

We can now implement Proposition~\ref{prop4.1} by defining
$h_{s_{1},\ldots,s_{k}}(y)$ to be the polynomial obtained from
$\phi_{s_1,\ldots,s_k}(y)$ by changing all subleading signs to $-$
as in \reff{eq4.1},
and letting $r(s_1,\ldots,s_k)$ be the unique positive root of
$h_{s_{1},\ldots,s_{k}}(y)$.
We then have $\rho(s_1,\ldots,s_k) \le r(s_1,\ldots,s_k)$.
Unfortunately, this bound is unsuitable for our present purposes,
as $r(s_1,\ldots,s_k)$ is not a monotonically decreasing
function of the path lengths $s_1,\ldots,s_k$
(see Table~\ref{table1new} and Remark 2 below).
We therefore throw away a bit more, by disregarding all sign cancellations
among the subleading terms of \reff{phiform}, and define
\be
\widetilde{h}_{s_{1},\ldots,s_{k}}(y)   \;=\;
      y^{(\sum{s_{i}}) - 1} \,-\,
      \sum_{m=2}^k \sum_{X \subseteq {[k] \choose {k-m}}}
           y^{s_X} (1 + y + \cdots + y^{m-2})
      \;.
 \label{hpoly}
\ee
[Thus, the coefficient of $y^j$ in \reff{hpoly}
 is in general {\em larger}\/ in magnitude than in \reff{phiform}.]
Let $\widetilde{r}(s_1,\ldots,s_k)$ be the unique positive root of
$\widetilde{h}_{s_{1},\ldots,s_{k}}(y)$.
Then it follows immediately from Proposition~\ref{prop4.1} that
\be
   \rho(s_1,\ldots,s_k)   \;\le\;   r(s_1,\ldots,s_k)
                          \;\le\;   \widetilde{r}(s_1,\ldots,s_k)
  \;,
\ee
or in other words:
\begin{prop}
\label{boundProp}
Every chromatic root $z$ of $\Theta_{s_{1},\ldots,s_{k}}$
lies in the disc $|z-1| \le \widetilde{r}(s_1,\ldots,s_k)$.
\end{prop}


We now analyze the behavior of the upper bound $\widetilde{r}(s_1,\ldots,s_k)$:

\begin{prop}
$\widetilde{r}(s_{1},\ldots ,s_{k})$ is symmetric in $s_1,\ldots,s_k$
and strictly decreasing in each $s_i$.
\label{dec}
\end{prop}

\proof
The symmetry is obvious.
To prove the decreasing property, fix $s_1,\ldots,s_k$
and set $r=\widetilde{r}(s_{1},s_{2},\ldots ,s_{k})$;
by symmetry, it suffices to show that $\widetilde{r}(s_{1}+1,s_{2},\ldots ,s_{k})<r$.
Now, it is clear from equation (\ref{hpoly})
that $r>1$. Note that from equation (\ref{hpoly})
we can rewrite $\widetilde{h}_{s_{1},s_{2},\ldots ,s_{k}}(y)$ in the form
\begin{eqnarray}
   \widetilde{h}_{s_{1},s_{2},\ldots,s_{k}}(y) & = & y^{(\sum_{i=1}^{k}
s_{i})-1} - y^{s_{1}}A(y) - B(y)
\label{hRewrite}
\end{eqnarray}
where $A$ and $B$ are polynomials in $y$ with nonnegative integer
coefficients that do not depend on $s_{1}$. Moreover, the degrees of
$y^{s_{1}}A(y)$ and $B(y)$ are less than $(\sum_{i=1}^{p}s_{i})-1$, and
$B(0)= 1$. It follows that
$r^{(\sum_{i=1}^{k}s_{i})}-r^{s_{1}}A(r) > 0$. Now from
(\ref{hRewrite}) we have
\begin{eqnarray}
\widetilde{h}_{s_{1}+1,s_{2},\ldots,s_{k}}(r) & = &
r^{\sum_{i=1}^{k} s_{i}} - r^{s_{1}+1}A(r) - B(r)   \nonumber \\
 & = & r \cdot \left( r^{(\sum_{i=1}^{k}s_{i}) -1} - r^{s_{1}} A(r)
\right) - B(r)  \nonumber \\
 & > & r^{(\sum_{i=1}^{k} s_{i})-1} - r^{s_{1}} A(r) - B(r)  \nonumber \\
 & = & \widetilde{h}_{s_{1},s_{2},\ldots,s_{k}}(r)  \nonumber \\
 & = & 0   \;.
\end{eqnarray}
It follows that $\widetilde{r}(s_{1}+1,s_{2},\ldots,s_{k}) < r =
\widetilde{r}(s_{1},s_{2},\ldots,s_{k})$, completing the proof. \qed

\medskip

What may be surprising is how well the roots of $h$ and $\widetilde{h}$
bound the roots of $f$ for small $k$ (see Table~\ref{table1new}).
In particular, they are considerably better than the bound
$\scrr(s_1,\ldots,s_k)$ from Theorem~\ref{thm3.2},
and they are good enough to prove:

\begin{table}[p]
\hspace*{-5mm}
\begin{tabular}{|l|c|c|c|c|} \hline
Path length sequence  & Actual value & \multicolumn{3}{|c|}{Upper bound} \\
\multicolumn{1}{|c|}{$(s_1,\ldots,s_k)$} & $\rho(s_1,\ldots,s_k)$ &
    $r(s_1,\ldots,s_k)$ & $\widetilde{r}(s_1,\ldots,s_k)$ &
    $\scrr(s_1,\ldots,s_k)$  \\
\hline\hline
(2, 2, 2)   &   1.5247025799   &   1.5905667405   &   1.5905667405   &   
   3.1478990357 \\
(2, 2, 3)   &   1.3247179572   &   1.4655712319   &   1.4655712319   &   
   2.8235871268 \\
\hline
(2, 2, 2, 2)   &   1.9635530390   &   2.0652388409   &   2.0959187459   &   
   3.6296581268 \\
(2, 2, 2, 3)   &   1.6180339887   &   1.8003794650   &   1.9038165409   &   
   3.3067093454 \\
\hline
(2, 2, 2, 2, 2)   &   2.3602010481   &   2.4788311017   &   2.5569445891
      &   4.0795956235 \\
(2, 2, 2, 2, 3)   &   1.9596554046   &   2.0481965587   &   2.3283569921
      &   3.7595287461 \\
(2, 2, 2, 2, 4)   &   1.9125157044   &   2.0726410424   &   2.2158195963
      &   3.6668970270 \\
(2, 2, 2, 2, 5)   &   2.0227195761   &   2.1137657905   &   2.1572723181
      &   3.6401168028 \\
(2, 2, 2, 2, 6)   &   1.9492237868   &   2.0928219450   &   2.1267590770
      &   3.6325613931 \\
\hline
(2, 2, 2, 2, 2, 2)   &   2.7305222731   &   2.8521866737   &   2.9891971006
      &   4.5063232460 \\
(2, 2, 2, 2, 2, 3)   &   2.3291754791   &   2.4702504048   &   2.7400794700
      &   4.1896653876 \\
(2, 2, 2, 2, 2, 4)   &   2.3208606055   &   2.4487347678   &   2.6342641478
      &   4.1075181051 \\
(2, 2, 2, 2, 3, 3)   &   2.0524815723   &   2.2641426827   &   2.5176585462
      &   3.8793014522 \\
\hline
(2, 2, 2, 2, 2, 2, 2)   &   3.0823336669   &   3.1959268744   &   3.4006086206
      &   4.9150761863 \\
(2, 2, 2, 2, 2, 2, 3)   &   2.6933092033   &   2.8543267466   &   3.1395749040
      &   4.6019501648 \\
(2, 2, 2, 2, 2, 2, 4)   &   2.7030241913   &   2.8316875864   &   3.0429807861
      &   4.5281826533 \\
(2, 2, 2, 2, 2, 3, 3)   &   2.3573224846   &   2.4527687226   &   2.8983449779
      &   4.2931001487 \\
\hline
(2, 2, 2, 2, 2, 2, 2, 2)   &   3.4201564280   &   3.5685068590   &   
   3.7959050193   &   5.3093300653 \\
(2, 2, 2, 2, 2, 2, 2, 3)   &   3.0446178232   &   3.2040479885   &   
   3.5278440533   &   4.9996840573 \\
(2, 2, 2, 2, 2, 2, 2, 4)   &   3.0625912820   &   3.2129169213   &   
   3.4402140830   &   4.9327412477 \\
(2, 2, 2, 2, 2, 2, 2, 5)   &   3.0953618332   &   3.1953189320   &   
   3.4125677445   &   4.9187003835 \\
(2, 2, 2, 2, 2, 2, 3, 3)   &   2.6885399588   &   2.8486049323   &   
   3.2745245420   &   4.6929626253 \\
\hline
(2, 2, 2, 2, 2, 2, 2, 2, 2)   &   3.7468849281   &   3.9272779941   &   
   4.1781887719   &   5.6915378807 \\
(2, 2, 2, 2, 2, 2, 2, 2, 3)   &   3.3836067543   &   3.5282506474   &   
   3.9060114610   &   5.3852446658 \\
(2, 2, 2, 2, 2, 2, 2, 2, 4)   &   3.4054981704   &   3.5867024115   &   
   3.8263498519   &   5.3239577745 \\
(2, 2, 2, 2, 2, 2, 2, 2, 5)   &   3.4292505541   &   3.5677746122   &   
   3.8040844502   &   5.3121036374 \\
(2, 2, 2, 2, 2, 2, 2, 2, 6)   &   3.4182415134   &   3.5704257784   &   
   3.7980747620   &   5.3098533475 \\
(2, 2, 2, 2, 2, 2, 2, 2, 7)   &   3.4200422197   &   3.5685857538   &   
   3.7964779130   &   5.3094286637 \\
(2, 2, 2, 2, 2, 2, 2, 2, 8)   &   3.4203605983   &   3.5684058522   &   
   3.7960560504   &   5.3093486377 \\
(2, 2, 2, 2, 2, 2, 2, 2, 9)   &   3.4200947731   &   3.5685249008   &   
   3.7959448158   &   5.3093335634 \\
(2, 2, 2, 2, 2, 2, 2, 2, 10)   &   3.4201605551   &   3.5685220773   &   
   3.7959155041   &   5.3093307241 \\
(2, 2, 2, 2, 2, 2, 2, 2, 11)   &   3.4201602535   &   3.5685079914   &   
   3.7959077815   &   5.3093301894 \\
(2, 2, 2, 2, 2, 2, 2, 2, 12)   &   3.4201547358   &   3.5685071412   &   
   3.7959057470   &   5.3093300886 \\
(2, 2, 2, 2, 2, 2, 2, 2, 13)   &   3.4201566935   &   3.5685072051   &   
   3.7959052110   &   5.3093300697 \\
(2, 2, 2, 2, 2, 2, 2, 3, 3)   &   3.0254986086   &   3.2079141314   &   
   3.6449248003   &   5.0809413850 \\
\hline
\end{tabular}
\caption{
   Values of $\rho(s_1,\ldots,s_k)$ and its upper bounds
   $r(s_1,\ldots,s_k) \le \widetilde{r}(s_1,\ldots,s_k)$.
   For comparison, the upper bound $\scrr(s_1,\ldots,s_k)$
   from Theorem~\ref{thm3.2} is also shown.
}
 \label{table1new}
\end{table}

\begin{thm}
  \label{2kbound}
For $3\leq k\leq 8$, we have $\rho(s_1,\ldots,s_k) \le \rho^{(2,k)}$,
with equality only when $s_1 = \cdots = s_k = 2$.
In other words, among all $k$-ary theta graphs,
the graph with a chromatic root that maximizes $|z-1|$
is the one with all path lengths equal to 2.
\end{thm}

\proof
By direct calculation (see Table~\ref{table1new})
we have $\widetilde{r}(2,2,3) < \rho(2,2,2)$,
$\widetilde{r}(2,2,2,3) < \rho(2,2,2,2)$ and
$\widetilde{r}(2,2,2,2,3) < \rho(2,2,2,2,2)$.
The result for $3 \le k \le 5$ then follows immediately
from Propositions~\ref{boundProp} and \ref{dec}.
For $k=6,7$, a bit more work is needed:
the calculations show that
$\rho(2,\ldots,2,2,3)$, $\widetilde{r}(2,\ldots,2,2,4)$ and
$\widetilde{r}(2,\ldots,2,3,3)$
are all bounded above by $\rho(2,\ldots,2,2,2)$,
so that the result again follows from
Propositions~\ref{boundProp} and \ref{dec}.
Finally, for $k=8$,
the calculations show that
$\rho(2,\ldots,2,2,3)$, $\rho(2,\ldots,2,2,4)$,
$\widetilde{r}(2,\ldots,2,2,5)$ and $\widetilde{r}(2,\ldots,2,3,3)$
are all bounded above by $\rho(2,\ldots,$ $2,2,2)$,
which is again sufficient.
\qed

\medskip
\noindent
{\bf Remarks.}
1.  This method of proof relies in an essential way on the fact that,
for $3 \le k \le 8$,
only finitely many of the upper bounds $\widetilde{r}(s_1,\ldots,s_k)$
are larger than the true value $\rho(2,\ldots,2) \equiv \rho^{(2,k)}$.
Unfortunately, this fails for $k=9$:
indeed, we have
$\lim\limits_{s_k \to \infty} \widetilde{r}(s_1,\ldots,s_{k-1},s_k) =
 \widetilde{r}(s_1,\ldots,s_{k-1})$
and in particular
$\lim\limits_{s_9 \to \infty} \widetilde{r}(2,\ldots,2,s_9) =
 \widetilde{r}^{(2,8)} \approx 3.7959050193 > 3.7468849281 \approx
 \rho^{(2,9)}$.
So a genuinely new method will be required to prove
Theorem~\ref{2kbound} for $k \ge 9$.

2.  Contrary to what one might expect,
neither the true value $\rho(s_1,\ldots,s_k)$
nor the upper bound $r(s_1,\ldots,s_k)$
is monotone decreasing in $s_1,\ldots,s_k$
(see Table~\ref{table1new}).
Nor does this arise merely from even-odd oscillations:
the weaker conjecture that
$\rho(s_1,\ldots,s_{k-1},s_k +2) \le \rho(s_1,\ldots,s_{k-1},s_k)$
for {\em even}\/ $s_1,\ldots,s_k \ge 2$ is also false,
as is the corresponding conjecture for $r$:
this is illustrated by the case
$(s_1,\ldots,s_5) = (2,2,2,2,4)$ among many others.

\section{The Chromatic Roots of $\Theta^{(2,k)} \simeq K_{2,k}$}  \label{K2k}

The results of the previous sections suggest that $\Theta^{(2,k)}$,
which is isomorphic to the complete bipartite graph $K_{2,k}$,
may very well contain a root that maximizes $|z-1|$ over all
$k$-ary theta graphs.
It is therefore of interest to study in detail the chromatic roots of $K_{2,k}$
and in particular their behavior as $k \to\infty$.
In this section we shall show that the bound
$\rho(s_1,\ldots,s_k) \le \scrr^{(2,k)} \approx [1+o(1)] \, k/\log k$
found in Theorem~\ref{thm3.2} and Proposition~\ref{prop_R2k}
is indeed asymptotically attained by $K_{2,k}$.
To do so, we shall need to carry out a rather in-depth study
of the roots of certain trinomials.

\begin{figure}[p]
\vspace*{0cm}
\begin{center}
\epsfxsize = 0.55\textwidth
\leavevmode\epsffile{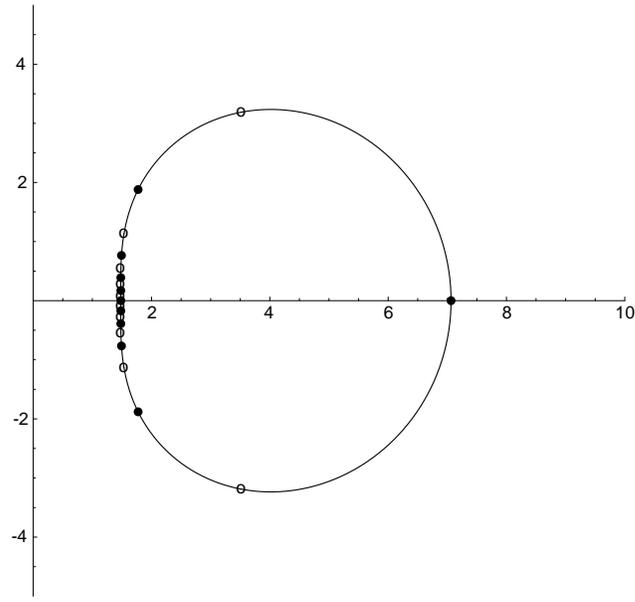} \qquad (a) \\
\vspace{1.5cm}
\epsfxsize = 0.55\textwidth
\leavevmode\epsffile{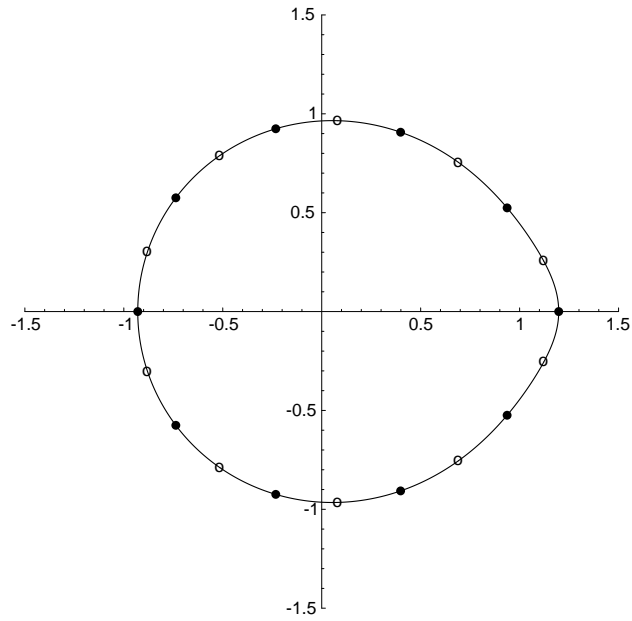} \qquad (b)
\end{center}
\vspace*{0.5cm}
\caption{
   Solution curve of the trinomial equation
   \reff{K2k_star1lam}/\reff{K2k_star3}
   for $|\lambda| = 1$, with $k=10$.
   Points correspond to $\lambda=+1$ ($\bullet$)
   and $\lambda=-1$ ($\circ$);
   the latter are the chromatic roots of the graph $K_{2,k}$.
   Plot (a) shows the complex $z$-plane;
   plot (b) shows the complex $\zeta$-plane.
}
 \label{fig2}
\end{figure}

{}From \reff{theta} we have
\be
   \pi(K_{2,k}, z)  \;=\;  z(z-1) [(z-2)^k + (z-1)^{k-1}]
   \;,
\ee
so that the chromatic roots of $K_{2,k}$ (aside from $z=0$ and $z=1$)
are given by the equation
\be
   (z-2)^k \,+\, (z-1)^{k-1}   \;=\;  0   \;.
 \label{K2k_star1}
\ee
By the Beraha--Kahane--Weiss theorem
\cite{BKW_75,BKW_78,Beraha_79,Beraha_80,Sokal_00b},
the solutions of \reff{K2k_star1} accumulate as $k \to\infty$
on the curve where $|z-2| = |z-1|$,
namely the vertical line $\Re z = 3/2$, and only there.
We were therefore somewhat surprised to compute numerically
the solutions of \reff{K2k_star1} and find that a few of the roots lie
{\em far to the right}\/ of the line $\Re z = 3/2$:
see e.g.\ the open-circled points in Figure~\ref{fig2}(a) for $k=10$.
In retrospect, however, one realizes that this behavior is
perfectly consistent with the Beraha--Kahane--Weiss theorem:
the limiting curve $\Re z = 3/2$ contains the point at infinity,
and roots can tend to infinity in the topology of the Riemann sphere
(which turns out to be the relevant sense)
in many different ways;
in particular, their real parts need not tend to 3/2.
In fact, as we shall see, the rightmost root has
\be
   z  \;=\;  {k \over \log k \,-\, \log\log k}
             \left[ 1 \,\pm\, {\pi i \over \log k}
                      \,+\, O\!\left( {\log\log k \over \log^2 k} \right)
             \right]
   \;,
 \label{K2k_star2}
\ee
so that its real and imaginary parts both {\em tend to infinity}\/
as $k \to\infty$,
with $|z|$ having exactly the magnitude $[1+o(1)] \, k/\log k$
predicted by Theorem~\ref{thm3.2} and Proposition~\ref{prop_R2k}.

The asymptotic behavior as $k \to\infty$ of the solutions to \reff{K2k_star1}
is surprisingly subtle,
and involves the Lambert $W$ function \cite{Corless_96}.
Here we shall give a brief treatment that leads as directly as possible
to the main result \reff{K2k_star2},
deferring a full analysis to a later paper \cite{Sokal_trinomials}
(see also \cite{dilcher,dilcher_unpub} for related work).

It is useful to study, in place of \reff{K2k_star1},
the more general equation
\be
   (z-2)^k \,-\, \lambda (z-1)^{k-1}   \;=\;  0
 \label{K2k_star1lam}
\ee
with any fixed $\lambda \neq 0$ ($\lambda \in \C$).
In particular, $\lambda = -1$ corresponds to the chromatic roots,
while the equation with $\lambda = +1$ has a positive real root
$z = 1 + \scrr^{(2,k)}$ where $\scrr^{(2,k)}$ is the bound
given in Theorem~\ref{thm3.2} [cf.\ \reff{eq3_Rc}].
See Figure~\ref{fig2}(a) for the curve of roots
corresponding to $|\lambda| = 1$, for the case $k=10$.

It is also convenient to make the fractional-linear change of variables
\be
   \zeta   \;=\;  {z-1 \over z-2}
   \;,
 \label{def_zeta}
\ee
so that the equation \reff{K2k_star1lam} becomes
\be
   \zeta^k - \zeta^{k-1} - \lambda   \;=\;  0
 \label{K2k_star3}
\ee
and the limiting curve $\Re z = 3/2$ becomes the unit circle $|\zeta| = 1$.
In particular, the point $z = \infty$ corresponds to $\zeta = 1$.
See Figure~\ref{fig2}(b) for the corresponding loci of roots in the
$\zeta$-plane, for the case $k=10$.

When $k \to\infty$ at fixed $\lambda \neq 0$, all the roots of
\reff{K2k_star3} tend to the unit circle $|\zeta| = 1$, as just noted;
but the {\em rate}\/ at which they do so is rather subtle,
most notably for those roots near $\zeta = 1$.
To investigate the large-$k$ asymptotics of the roots of \reff{K2k_star3}
near $\zeta = 1$, let us begin by writing $\zeta = e^{w/k}$
with $|w/k| \ll 1$ and inserting this in \reff{K2k_star3}:  we have
\be
   w e^w  \;=\; k \lambda
                \,  \left[ 1 \,+\, O\Bigl({w \over k}\Bigr) \right]
   \;.
\ee
This suggests that to a first approximation we will have
\be
   w  \;=\;  W(k\lambda)
 \label{def_w}
\ee
where $W(z)$ is the Lambert $W$ function \cite{Corless_96},
i.e.\ the (multivalued) inverse function to $w \mapsto w e^w$,
and distinct branches of $W(z)$ will correspond to distinct roots
of the trinomial.
The large-$z$ behavior of $W(z)$ is rather complicated \cite{Corless_96},
but the first two terms of its expansion are
\be
   W(z) \;=\; \log z \,-\, \log \log z \,+\,
              O\!\left( {\log \log z \over \log z} \right)
   \;,
 \label{W_large_z}
\ee
so that indeed we will have $|w/k| \approx |\log k/k| \ll 1$
for large $k$ (provided that $\arg\lambda$ is not taken too large).

Let us henceforth {\em define}\/ $w$ by \reff{def_w},
using any branch of the Lambert $W$ function we please,
and write
\be
   \zeta   \;=\;   \exp\!\left[ {w \over k} (1+\xi) \right]
   \;,
 \label{def_xi}
\ee
where we expect $|\xi| \ll 1$ for large $k$.
Inserting this into \reff{K2k_star3},
we obtain the {\em fundamental equation}\/
\be
   {1 \,-\, e^{-\tau(1+\xi)}  \over  \tau}  \;=\;  e^{-w\xi}
   \;,
 \label{fund_eqn_w}
\ee
where we have set
\be
   \tau   \;=\;  {w \over k}
   \;.
 \label{def_tau}
\ee
Our strategy will now be to investigate the equation \reff{fund_eqn_w},
with $w$ and $\tau$ considered as independent complex parameters,
and seek a convergent power-series solution
\be
   \xi(w,\tau)  \;=\;  \sum\limits_{l=1}^\infty f_l(w) \, \tau^l
 \label{series_xi_w}
\ee
valid for $|\tau| < A(w)$.  Indeed, the implicit function theorem
guarantees that such a solution exists for small enough $|\tau|$,
provided that $w \neq -1$.
In view of the special role played by $w=-1$,
it is convenient to make the further change of variables
\be
   v  \;=\;  {1 \over 1+w}
   \;,
 \label{def_v}
\ee
so that the fundamental equation \reff{fund_eqn_w} becomes
\be
   {1 \,-\, e^{-\tau(1+\xi)}  \over  \tau}  \;=\;  e^{(1-1/v)\xi}
 \label{fund_eqn_v}
\ee
and the power-series solution will be
\be
   \xi(v,\tau)  \;=\;  \sum\limits_{l=1}^\infty g_l(v) \, \tau^l
   \;.
 \label{series_xi_v}
\ee

Starting from \reff{fund_eqn_v}, let us take the logarithm of both sides,
subtract $\xi$, and multiply by $-v$:  we get
\be
   \xi   \;=\;
   -v \left[ \log\!\left( {1 \,-\, e^{-\tau(1+\xi)}  \over  \tau(1+\xi)}
                   \right)
             \,+\, \log(1+\xi) \,-\, \xi
      \right]
   \;.
 \label{eq5.xxx}
\ee
We intend to prove that, whenever $|\tau| \le A$ and $|v| \le B$,
the equation \reff{eq5.xxx} has a unique solution $\xi = \xi(v,\tau)$
lying in the disc $|\xi| \le R$,
provided that $A,B,R$ satisfy suitable inequalities.
To prove this, we shall apply Rouch\'e's theorem to \reff{eq5.xxx},
taking the left-hand side as the ``large'' function $f(\xi)$ and the
right-hand side as the ``small'' function $g(\xi)$.
Trivially $f(\xi) = \xi$ is analytic in the entire $\xi$-plane,
and its only zero is a simple zero at $\xi=0$.
The right-hand side $g(\xi)$ is analytic in the disc $|\xi| < R$
and continuous in $|\xi| \le R$, provided that
$A < 2\pi$ and $R < \min(1, 2\pi/A - 1)$.
Moreover, we can bound $g(\xi)$ in this disc by observing that
\begin{subeqnarray}
   \log\!\left( {1 - e^{-z} \over z} \right)
   & = &   - {z \over 2}  \,+\,  \log\!\left( {\sinh(z/2) \over z/2} \right)
      \\[2mm]
   & = &   - {z \over 2}  \,+\,
           \sum\limits_{n=1}^\infty {(-1)^{n-1} \zeta(2n) \over n (2\pi)^{2n}}
                \, z^{2n}
   \quad\hbox{for } |z| < 2\pi
 \label{thesinh}
\end{subeqnarray}
and hence
\be
   \left| \log\!\left( {1 - e^{-z} \over z} \right) \,+\, {z \over 2} \right|
   \;\le\;
   \sum\limits_{n=1}^\infty {\zeta(2n) \over n (2\pi)^{2n}} \, |z|^{2n}
   \;=\;
   \log\!\left( {|z|/2 \over \sin(|z|/2)} \right)
   \quad\hbox{for } |z| < 2\pi
   \;,
\ee
and similarly
\be
   \log(1+\xi) \,-\, \xi   \;=\;
      \sum\limits_{n=2}^\infty {(-1)^{n-1} \over n} \, \xi^n
   \quad\hbox{for } |\xi| < 1
\ee
and hence
\be
   |\log(1+\xi) \,-\, \xi|
   \;\le\;
   \sum\limits_{n=2}^\infty {|\xi|^n \over n}
   \;=\;
   -\log(1-|\xi|) \,-\, |\xi|
   \quad\hbox{for } |\xi| < 1
   \;.
\ee
Therefore, for $|\tau| \le A < 2\pi$, $|v| \le B$ and
$|\xi| \le R < \min(1, 2\pi/A - 1)$, we have
\be
   |g(\xi)|   \;\le\;
   B \left[  {A(1+R) \over 2}
             \,+\,
             \log\!\left( {A(1+R)/2 \over \sin[A(1+R)/2]} \right)
             \,-\, \log(1-R) \,-\, R
     \right]
   \;.
 \label{bound_g}
\ee
Rouch\'e's theorem applies provided that $|f(\xi)| > |g(\xi)|$
everywhere on the circle $|\xi| = R$;
a sufficient condition for this is thus
\be
   B \left[  {A(1+R) \over 2}
             \,+\,
             \log\!\left( {A(1+R)/2 \over \sin[A(1+R)/2]} \right)
             \,-\, \log(1-R) \,-\, R
     \right]
   \;<\;  R
   \;.
 \label{Rouche_condition}
\ee

It is easy to see that whenever $A$ and $B$ are sufficiently small,
there exists an $R > 0$ satisfying \reff{Rouche_condition}.\footnote{
   Much more can be said, but we defer a complete analysis
   to a later paper \cite{Sokal_trinomials}.
}
So let $(A,B,R)$ be a triplet satisfying \reff{Rouche_condition}.
Then, whenever $|\tau| \le A$ and $|v| \le B$,
the equation \reff{eq5.xxx} has a unique solution $\xi = \xi(v,\tau)$
lying in the disc $|\xi| \le R$;
moreover, the implicit function theorem guarantees that
this solution is an analytic function of $\tau$ and $v$ in the open polydisc
$D_{A,B} \equiv \{ (\tau,v) \colon\;  |\tau| < A ,\, |v| < B \}$
(a simple zero moves analytically under small analytic perturbations).
In particular, $\xi(v,\tau)$ is given in $D_{A,B}$ by
an absolutely convergent Taylor series
\be
   \xi(v,\tau)  \;=\;
   \sum\limits_{l=1}^\infty \sum\limits_{m=1}^\infty g_{lm} \tau^l v^m
   \;=\;  \sum\limits_{l=1}^\infty g_l(v) \, \tau^l
   \;.
 \label{series2_xi_v}
\ee
The first few $g_l(v)$ are easily computed by expansion of \reff{fund_eqn_v}:
\begin{subeqnarray}
   g_1(v)  & = &    v/2   \slabel{g1_v} \\[2mm]
   g_2(v)  & = &    (-v + 6v^2 + 3v^3)/24   \\[2mm]
   g_3(v)  & = &    (-3v^2 + 5v^3 + 7v^4 + 3v^5)/48   \\[2mm]
           & \vdots &   \nonumber
\end{subeqnarray}

It easily follows from the convergence of \reff{series2_xi_v}
that for $|\tau| \le A' < A$ and $|v| \le B' < B$ we have
\be
   |\xi(v,\tau)|  \;\le\;  C \, |\tau| \, |v|
                  \;=\;  {C \over k}  \left| {w \over 1+w} \right|
\ee
for a suitable constant $C < \infty$.
So let $\theta$ be any fixed real number,
and let $\lambda = e^{i\theta}$ (considered as belonging to
the Riemann surface of the logarithm function).
If we take $k$ large enough (how large depends on $\theta$),
then $w \equiv W(k e^{i\theta})$ will satisfy [by \reff{W_large_z}]
\be
   w  \;=\;  \log k \,+\, i\theta \,-\, \log\log k
             \,+\, O\!\left( {\log\log k \over \log k} \right)
\ee
and hence in particular $|\tau| = |w/k| \le A'$
and $|v| = |1+w|^{-1} \le B'$.
It follows that the trinomial \reff{K2k_star3} has a solution $\zeta$
satisfying
\begin{subeqnarray}
   \zeta  & = &   \exp\!\left[ {w \over k} \,+\,
                               O\!\left( {w \over k^2} \right)
                        \right]
   \\[2mm]
   & = &   \exp\!\left[ {\log k \over k} \,+\, {i\theta \over k}
                        \,-\, {\log\log k \over k}
                        \,+\, O\!\left( {\log\log k \over k \log k} \right)
                        \right]
   \\[2mm]
   & = &  1 \,+\, {\log k \over k} \,+\, {i\theta \over k}
            \,-\, {\log\log k \over k}
            \,+\, O\!\left( {\log\log k \over k \log k} \right)
   \;.
\end{subeqnarray}
Transforming back to the variable $z = (2\zeta-1)/(\zeta-1)$,
we find
\begin{subeqnarray}
   z   & = &  {k \over  \log k \,+\, i\theta \,-\, \log\log k
                           \,+\, O\!\left( {\log\log k \over \log k} \right)
              }
   \\[2mm]
   & = &  {k \over  \log k \,-\, \log\log k}
          \, \left[ 1 \,-\, {i\theta \over \log k}
                      \,+\, O\!\left( {\log\log k \over \log^2 k} \right)
             \right]
          \;.
\end{subeqnarray}
We have therefore proven:

\begin{theorem}
   \label{thm5.1}
Fix $\theta \in \R$.  Then, for all sufficiently large $k$,
the equation $(z-2)^k - e^{i\theta} (z-1)^{k-1} = 0$ has a solution
\be
   z   \;=\;
   {k \over  \log k \,-\, \log\log k}
   \, \left[ 1 \,-\, {i\theta \over \log k}
               \,+\, O\!\left( {\log\log k \over \log^2 k} \right)
      \right]
   \;.
\ee
\end{theorem}

\begin{corollary}
   \label{cor5.2}
For all sufficiently large $k$, the graph $\Theta^{(2,k)} \simeq K_{2,k}$
has a pair of chromatic roots
\be
   z  \;=\;  {k \over \log k \,-\, \log\log k}
             \left[ 1 \,\pm\, {\pi i \over \log k}
                      \,+\, O\!\left( {\log\log k \over \log^2 k} \right)
             \right]
   \;.
\ee
\end{corollary}

\section{Concluding Remarks}  \label{sec_conclusion}

A $k$-ary generalized theta graph is a series-parallel (hence planar)
graph of maximum degree $k$ and corank $k-1$
(except for the trivial case $\Theta_{s_1}$ with $s_1 \ge 2$).
Our main result, Theorem~\ref{thm1.3}, therefore naturally suggests
extensions in two different directions:
\begin{itemize}
   \item[(a)]  a sublinear bound in terms of maximum degree
for larger classes of series-parallel (or perhaps even planar) graphs, and
   \item[(b)]  a sublinear bound in terms of corank for larger classes
of graphs, and possibly even for arbitrary graphs.
\end{itemize}
We discuss these conjectures in turn:

\subsection{Bounds in terms of maximum degree}

Theorem~\ref{thm1.1} provides a linear bound in terms of maximum degree
for the chromatic roots of arbitrary graphs.
For general graphs, such a bound is best possible
except for the numerical value of the prefactor,
since the complete graph $K_{k+1}$ has a chromatic root at $z=k$.
On the other hand, suitably restricted subclasses of graphs
might well satisfy a sublinear bound,
as we have shown in Theorem~\ref{thm1.3} for generalized theta graphs,
whose chromatic roots are bounded by $[1 + o(1)] \, k/\log k$.
We conjecture that the techniques of Section~\ref{bound}
can be extended to handle arbitrary {\em series-parallel}\/ graphs:

\begin{conjecture}
There exists a universal constant $C < \infty$
such that the chromatic roots of any series-parallel graph
of maximum degree $k$ lie in the disc $|z-1| \le C k/\log k$.
\end{conjecture}

\noindent
More strongly, we conjecture that the chromatic roots lie in the
same disc $|z-1| \le \scrr^{(2,k)}$ that we have established
in Theorem~\ref{thm3.2} for generalized theta graphs.
Indeed, it is quite possible that among series-parallel graphs
the worst case is always the generalized theta graph
$\Theta^{(2,k)} \simeq K_{2,k}$,
so that the roots lie in the disc $|z-1| \le \rho^{(2,k)}$;
this generalizes the conjectured extension of Theorem~\ref{2kbound}
to $k \ge 9$.

Series-parallel graphs are a subset of planar graphs,
so it is conceivable that a sublinear bound in terms of maximum degree
holds even for all {\em planar}\/ graphs.
But we have no idea how to prove such a result,
nor do we have any compelling reason to believe it is true.

Another direction in which Theorems~\ref{thm1.1} and \ref{thm1.3}
could be extended is by finding a criterion {\em weaker}\/
than bounded maximum degree under which the chromatic roots
could be shown to be bounded (whether linearly or sublinearly).
Indeed, already in \cite{Sokal_00a} it was shown ``maximum degree''
in Theorem~\ref{thm1.1} can be replaced by ``second-largest degree'',
provided that the bound $7.963907\,k$
is replaced by $7.963907\,k \,+\, 1$.\footnote{
   This result {\em cannot}\/ be extended further to ``third-largest degree'':
   for as soon as $G$ has {\em two}\/ vertices of large degree,
   the chromatic roots can become unbounded.
   Indeed, the chromatic roots of the generalized theta graphs $\Theta^{(s,k)}$,
   when $s$ and $k$ are both allowed to vary without bound,
   are dense in the entire complex plane except possibly for the disc
   $|z-1| < 1$ \cite{Sokal_00b}.
}
And it was conjectured there (inspired by \cite{Shrock_99a})
that ``second-largest degree'' can be further weakened to ``maxmaxflow'',
defined as
\be
   \Lambda(G)   \;=\;   \max\limits_{x \neq y}  \lambda(x,y)
\ee
where
\begin{subeqnarray}
   \lambda(x,y)
      & = &  \hbox{max \# of edge-disjoint paths from $x$ to $y$} \\
      & = &  \hbox{min \# of edges separating $x$ from $y$}
\end{subeqnarray}
(Clearly $\lambda(x,y) \le \min[\deg(x), \deg(y)]$
and hence $\Lambda(G) \le $ second-largest degree of $G$.)
In other words, it was conjectured that

\begin{conjecture}[\cite{Shrock_99a,Sokal_00a}]
There exist universal constants $C(k) < \infty$
such that the chromatic roots of any graph $G$ with $\Lambda(G) = k$
lie in the disc $|z-1| \le C(k)$.
\end{conjecture}

\noindent
Indeed, it is natural to expect that $C(k)$ can be taken to be
{\em linear}\/ in $k$.
We do not have, at present, any good idea how to prove this conjecture
for arbitrary graphs, but we suspect that the methods of Section~\ref{bound}
can be extended to prove it for {\em series-parallel}\/ graphs.

\subsection{Bounds in terms of corank}

Recalling the definition $\rho(G) = \max\{ |z-1| \colon\; \pi(G,z) = 0 \}$,
let us now define the numbers
\be
   \rho_{k}  \;=\;  \max\{ \rho(G) \colon\; \mbox{graphs $G$ of corank $k$} \}
   \;.
 \label{eq6b.1}
\ee
Obviously $\rho_0 = \rho_1 = 1$.
Note also that $\rho_{k+1} \geq \rho_{k}$:
for if $G$ is any graph of corank $k$,
then the disjoint union of $G$ and a cycle $C_n$ has corank $k+1$
(and its chromatic roots are those of $G$ and the cycle).
Note, finally, that in the definition of $\rho_{k}$
we can restrict attention to 2-connected graphs:
for we can separate $G$ into its 2-connected components
and then glue these components back together along a common edge
to form a graph $G'$ of the same corank as $G$;
and by the $l=0,1,2$ cases of \reff{union},
the chromatic roots of $G'$ are exactly those of $G$
(except that the multiplicities of the roots at 0 and 1 may be reduced).

Theorem~\ref{thm1.2} asserts that $\rho_k \le k$ for $k \ge 1$.
In the other direction, Theorem~\ref{thm1.3} and Corollary~\ref{cor5.2}
together show that $\rho_k \ge \rho^{(2,k+1)} = [1 + o(1)] \, k/\log k$.
It is clearly of interest to know whether the asymptotic growth of $\rho_k$
is linear or sublinear.

Let us begin by examining small values of $k$.
As noted above, $\rho_0 = \rho_1 = 1$.
The only 2-connected graphs of corank 2 are 3-ary theta graphs,
so Theorem~\ref{2kbound} implies that $\rho_{2} = \rho(2,2,2) \approx 1.5247$.
We initially conjectured that an analogous result might hold for all $k$,
i.e.\ that the corank-$k$ graph with the largest value of $\rho(G)$
would be the $(k+1)$-ary theta graph $\Theta^{(2,k+1)}$,
so that $\rho_k$ would equal $\rho(2,\ldots,2) = \rho^{(2,k+1)}$.
Sadly, this is not the case, as $\rho_3 \geq 2$
(since $z=3$ is a chromatic root of $K_{4}$),
while $\rho(2,2,2,2) \approx 1.9636$ (see Table~\ref{table1new}).
However, this is the only counterexample we have found,
so we pose the following conjecture:

\begin{conj}
Let $G$ be a a graph of corank $k \ge 1$.
Then, if $G \neq K_4$, we have $\rho(G) \le \rho^{(2,k+1)}$.
In particular, for $k \ge 4$ we have $\rho_k = \rho^{(2,k+1)}$.
 \label{conj6.1}
\end{conj}

\noindent
In particular, we expect:

\begin{conj}[a corollary of Conjecture~\ref{conj6.1}]
$\rho_{k} = [1 + o(1)] \, k / \log k$ as $k \to\infty$.
\end{conj}

\noindent
Note that asymptotically, complete graphs will lag far behind
generalized theta graphs, as $\rho(K_n) = n-2$,
and the corank of $K_{n}$ is $(n^2-3n+2)/2$;
the corresponding generalized theta graph with the same corank,
$\Theta^{(2,1+(n^2-3n+2)/2)} \simeq K_{2,1+(n^2-3n+2)/2}$,
has, from Corollary~\ref{cor5.2},
a chromatic root $z$ such that $|z-1|$ is approximately $n^2/(4\log n)$,
which is much larger than $n-2$.

\section*{Acknowledgments}

We wish to thank Karl Dilcher for many helpful conversations
concerning the roots of trinomials,
and for sharing with us his unpublished notes \cite{dilcher_unpub}.
One of us (A.D.S.)\ also wishes to thank Rob Corless and David Jeffrey
for extremely helpful correspondence concerning the Lambert $W$ function;
in particular, Rob's computation of the first few terms of a series
closely related to \reff{series_xi_w}
was an essential stimulus for our realizing that such expansions
might be {\em convergent}\/.

This research was supported in part by NSF grant PHY--9900769 (A.D.S.)\ 
and operating grants from NSERC (J.B., C.H.\ and D.G.W.).
It was completed while one of the authors (A.D.S.)\ was a
Visiting Fellow at All Souls College, Oxford,
where his work was supported in part by EPSRC grant GR/M 71626
and aided by the warm hospitality of John Cardy and the
Department of Theoretical Physics.

\end{document}